\documentclass{amsart}
\usepackage{amssymb}
\usepackage{amsmath}
\usepackage{amsfonts}

\newtheorem{theorem}{Theorem}[section]
\theoremstyle{plain}

\newtheorem{corollary}[theorem]{Corollary}

\newtheorem{example}[theorem]{Example}

\newtheorem{lemma}[theorem]{Lemma}

\numberwithin{equation}{section}

\begin{document}
\title{Chief factors of Lie Algebras}
\author{David A. TOWERS}
\address{Lancaster University\\
Department of Mathematics and Statistics \\
LA$1$ $4$YF Lancaster\\
ENGLAND}
\email{d.towers@lancaster.ac.uk}

\thanks{$2010$ \textit{Mathematics Subject Classification.} $17$B$05$, $17$B$%
20$, $17$B$30$, $17$B$50$. }

\begin{abstract} In group theory the chief factors allow a group to be studied by its representation theory on particularly natural irreducible modules. It is to be expected, therefore, that they will play an important role in the study of Lie algebras. In this article we survey a few of their properties.
\end{abstract}

\maketitle

\section{Introduction}
Throughout $L$ will denote a finite-dimensional Lie algebra over a field $F$.  We call a subalgebra $I$ a {\em subideal} of a Lie algebra $L$ if there is a chain of subalgebras
\[ I=I_0 < I_1 < \ldots < I_n=L,
\] where $I_j$ is an ideal of $I_{j+1}$ for each $0 \leq j \leq n-1$.
\par

Put $L^1=L, L^{k+1}=[L^k,L]$ for $k \geq 1$. These are the terms of the lower central series for $L$. We say that $L$ has {\em nilpotency class $n$} if $L^n \neq 0$ but $L^{n+1}=0$. Let $U$ be a subalgebra of $L$. If $F$ has characteristic $p>0$ we call $U$ {\em nilregular} if the nilradical of $U$, $N(U)$, has nilpotency class less than $p-1$. If $F$ has characteristic zero we regard every subalgebra of $L$ as being nilregular. We say that $U$ is {\em characteristic in $L$} if it is invariant under all derivations of $L$. Nilregular ideals of $L$ have the property that their nilradicals are characteristic in $L$. Details of the results in this section can be found in \cite{gennil}.

\begin{theorem}\label{t:nil} \begin{itemize}
\item[(i)] If $I$ is a nilregular ideal of $L$ then $N(I) \subseteq N(L)$.
\item[(ii)] If $I$ is a nilregular subideal of $L$ and every subideal of $L$ containing $I$ is nilregular, then $N(I) \subseteq N(L)$.
\end{itemize}
\end{theorem}
\par

This result was proved by Schenkman (\cite{schenk}) for fields of characteristic zero; in characteristic $p$ it follows from a more recent result of Maksimenko (\cite{mak}). Similarly, we will call the subalgebra $U$ {\em solregular} if the underlying field $F$ has characteristic zero, or if it has characteristic $p$ and the (solvable) radical of $U$, $R(U)$, has derived length less than $log_2 p$. Then we have the following corresponding theorem, which uses a result of Petravchuk (\cite{pet}).

\begin{theorem}\label{t:solv}  \begin{itemize}
\item[(i)] If $I$ is a solregular ideal of $L$ then $R(I) \subseteq R(L)$.
\item[(ii)] If $I$ is a solregular subideal of $L$ and every subideal of $L$ containing $I$ is solregular, then $R(I) \subseteq R(L)$.
\end{itemize}
\end{theorem}
\medskip

These enable us to determine what the minimal ideals of $L$ look like.

\begin{theorem}\label{t:minp} Let $L$ be a Lie algebra over a field $F$, and let $I$ be a minimal non-abelian ideal of $L$. Then either 
\begin{itemize}
\item[(i)] $I$ is simple or 
\item[(ii)] $F$ has characteristic $p$, $N(I)$  has nilpotency class greater than or equal to $p-1$, and $R(I)$ has derived length greater than  or equal to $log_2 p$.
\end{itemize}
\end{theorem}

As a result of the above we will call the subalgebra $U$ {\em regular} if it is either nilregular or solregular; otherwise we say that it is {\em irregular}. 
Then we have the following corollary.

\begin{corollary}\label{c:minp} Let $L$ be a Lie algebra over a field $F$. Then every minimal ideal of $L$ is abelian, simple or irregular.
\end{corollary}

Block's Theorem on differentiably simple rings (see \cite{block}) describes the irregular minimal ideals as follows.

\begin{theorem}\label{t:block} Let $L$ be a Lie algebra over a field of characteristic $p>0$ and let $I$ be an irregular minimal ideal of $L$. Then $I \cong  S\otimes {\mathcal O}_n$, where $S$ is simple and ${\mathcal O}_n$ is the truncated polynomial algebra in $n$ indeterminates. Moreover, $N(I)$ has nilpotency class $p-1$ and $R(I)$ has derived length $\lceil log_2 p \rceil$.
\end{theorem}

\section{Primitive Lie algebras}
Next we  introduce the concept of a primitive Lie algebra. Details of the results in this section can be found in \cite{prim}. A word of warning - this terminology has been used for a different concept elsewhere. If $U$ is a subalgebra of $L$ we define $U_L$, the {\em core} (with respect to $L$) of $U$ to be the largest ideal of $L$ contained in $U$. We say that $U$ is {\em core-free} in $L$ if $U_L = 0$. We shall call $L$ {\em primitive} if it has a core-free maximal subalgebra. The {\em centraliser} of $U$ in $L$ is $C_L(U)= \{x \in L : [x,U]=0 \}$.
\par

There are three types of primitive Lie algebra: 
\begin{itemize}
\item[1.] {\em primitive of type $1$} if it has a unique minimal ideal that is abelian;
\item[2.] {\em primitive of type $2$} if it has a unique minimal ideal that is  non-abelian; and
\item[3.] {\em primitive of type $3$} if it has precisely two distinct minimal ideals each of which is  non-abelian.
\end{itemize}

Of course, primitive Lie algebras of types $2$ and $3$ are semisimple, and those of types $1$ and $2$ are monolithic. (A Lie algebra $L$ is called {\em monolithic} if it has a unique minimal ideal $W$, the {\em monolith} of $L$.)

\begin{example} Examples of each type are easy to find.
\begin{itemize}
\item[1.] Clearly every primitive solvable Lie algebra is of type $1$.
\item[2.] Every simple Lie algebra is primitive of type $2$.
\item[3.] If $S$ is a simple Lie algebra then $L=S \oplus S$ is primitive of type $3$ with core-free maximal subalgebra $D = \{ s+s : s \in S \}$, the diagonal subalgebra of $L$.
\end{itemize}  
\end{example}

Let $M$ be a maximal subalgebra of $L$. Then $M/M_L$ is a core-free maximal subalgebra of $L/M_L$. We say that $M$ is
\begin{itemize}
\item[1.] {\em a maximal subalgebra of type $1$} if $L/M_L$ is primitive of type $1$;
\item[2.] {\em a maximal subalgebra of type $2$} if  $L/M_L$ is primitive of type $2$; and
\item[3.] {\em a maximal subalgebra of type $3$} if  $L/M_L$ is primitive of type $3$.
\end{itemize}

We say that an ideal $A$ is {\em complemented} in $L$ if there is a subalgebra $U$ of $L$ such that $L=A+U$ and $A \cap U=0$. For primitive solvable Lie algebras we have the following analogue of Galois' Theorem for groups.

\begin{theorem}\label{t:galois}
\begin{itemize}
\item[1.] If $L$ is a solvable primitive Lie algebra then all core-free maximal subalgebras are conjugate.
\item[2.] If $A$ is a self-centralising minimal ideal of a solvable Lie algebra $L$, then $L$ is primitive, $A$ is complemented in $L$, and all complements are conjugate.
\end{itemize}
\end{theorem}
\medskip

 The {\em Frattini ideal} of $L$, $\phi(L)$, is the  core of intersection of the maximal subalgebras of $L$. We say that $L$ is {\em $\phi$-free} if $\phi(L)=0$. Then we have the following characterisation of primitive Lie algebras of type $1$.

\begin{theorem}\label{t:type1} Let $L$ be a Lie algebra over a field $F$.
\begin{itemize}
\item[1.] $L$ is primitive of type $1$ if and only if $L$ is monolithic, with abelian monolith $W$, and $\phi$-free.
\item[2.] If $F$ has characteristic zero, then $L$ is primitive of type $1$ if and only if $L=W \ltimes (C \oplus S)$ (semi-direct sum), where $W$ is the abelian monolith of $L$, $C$ is an abelian subalgebra of $L$, every element of which acts semisimply on $W$, and $S$ is a Levi subalgebra of $L$.
\item[3.] If $L$ is solvable, then $L$ is primitive if and only if it has a self-centralising minimal ideal $A$. 
\end{itemize}
\end{theorem}
\medskip

For type $2$ we have

\begin{theorem}\label{t:type2}
\begin{itemize}
\item[1.] $L$ is primitive of type $2$ if and only if  $L \cong U + (S \otimes {\mathcal O}_n)$,
where $S \otimes {\mathcal O}_n$ is an ideal of $L$ and $S$ is simple. 
\item[2.] If $F$ has characteristic zero, then $L$ is primitive of type $2$ if and only if $L$ is simple.
\item[3.] $L$ is primitive of type $2$ if and only if there is a primitive Lie algebra $X$ of type $3$ such that $L \cong X/B$ for a minimal ideal $B$ of $L$. 
\end{itemize}
\end{theorem}
\medskip

For type $3$ we have

\begin{theorem}\label{c:type3}
\begin{itemize}
\item[1.] $L$ is primitive of type $3$ if and only if $L$ has two distinct minimal ideals $B_1$ and $B_2$ with a common complement and such that the factor algebras $L/B_i$ are primitive of type $2$ for $i=1,2$. Moreover, $B_1$ and $B_2$ are both isomorphic to $S \otimes {\mathcal O}_n$, where $S$ is simple. 
\item[2.] If $F$ has characteristic zero, then $L$ is primitive of type $3$ if and only if $L = S \oplus S$, where $S$ is simple.
\end{itemize}
\end{theorem}

\section{Chief factors}
 The factor algebra $A/B$ is called a {\em chief factor} of $L$ if $B$ is an ideal of $L$ and $A/B$ is a minimal ideal of $L/B$. So chief factors are as described in Corollary \ref{c:minp} and Theorem \ref{t:block}. We can identify different types of chief factor; details for this section can be found in \cite{prim}. A chief factor $A/B$ is called {\em Frattini }if $A/B\subseteq \phi
\left( L/B\right) .$ This concept was first introduced in \cite{prefrat}.
\par

If there is a subalgebra, $M$ such that $L=A+M$ and $B\subseteq A\cap M,$ we
say that $A/B\ $is$\ $a \textit{supplemented} chief factor of $L,$ and that $%
M$ is a \textit{supplement} of $A/B\ $in $L.$ Also, if $A/B$ is a
non-Frattini chief factor of $L$, then $A/B$ is \textit{supplemented} by a
maximal subalgebra $M$ of $L.$
\par

If $A/B\ $is$\ $a chief factor of $L$ supplemented by a subalgebra $M$ of $%
L, $ and \linebreak $A\cap M=B$ then we say that $A/B$ is \textit{%
complemented} chief factor of $L,$ and $M$ is a \textit{complement} of $A/B$
in $L$. When $L$ is \textit{solvable}, it is easy to see that a chief factor
is Frattini  if and only if it is not complemented. Then we have the following generalisation of the Jordan-H\"{o}lder Theorem.

\begin{theorem}\label{t:jordan} Let
\[ 0 < A_1 < \ldots < A_n = L  \hspace{1in} (1)
\]
\[ 0 < B_1 < \ldots < B_n = L \hspace{1in} (2)
\]
be chief series for the Lie algebra $L$. Then there is a bijection between the chief factors of these two series such that corresponding factors are isomorphic as $L$-modules and such that the Frattini chief factors in the two series correspond. 
\end{theorem}
\medskip
 The number of Frattini chief factors or of chief factors which are complemented by a maximal subalgebra of a finite-dimensional Lie algebra $L$ is the same in every chief series for $L$. However, this is not the case for the number of chief factors which are simply complemented in $L$; in \cite{tz} we determine the possible variation in that number.
\par

Note that if $L$ is a primitive Lie algebra of type $3$, its two minimal ideals are not $L$-isomorphic, so we introduce the following concept. We say that two chief factors of $L$ are {\em $L$-connected} if either they are $L$-isomorphic, or there exists an epimorphic image $\overline{L}$ of $L$ which is primitive of type $3$ and whose minimal ideals are $L$-isomorphic, respectively, to the given factors. (It is clear that, if two chief factors of $L$ are $L$-connected and are not $L$-isomorphic, then they are nonabelian and there is a single epimorphic image of $L$ which is primitive of type $3$ and which connects them.) Then, as we would hope,

\begin{theorem}\label{t:Lcon} The relation `is $L$-connected to' is an equivalence relation on the set of chief factors.
\end{theorem}
\medskip

Let $A/B$ be a supplemented chief factor of $L$ and put ${\mathcal J} =  \{M_L : M$ is a maximal subalgebra of $L$ supplementing a chief factor $L$-connected to $A/B \}$. Let $R= \cap\{N : N \in {\mathcal J} \}$ and $C=A+C_L(A/B)$. Then we call $C/R$ the {\em crown} of $L$ associated with $A/B$. This object gives much information about the supplemented chief factors of $L$.

\begin{theorem}\label{t:crown} Let $C/R$ be the crown associated with the supplemented chief factor $A/B$ of $L$. Then $C/R= Soc(L/R)$. Furthermore
\begin{itemize}
\item[(i)] every minimal ideal of $L/R$ is a supplemented chief factor of $L$ which is $L$-connected to $A/B$, and
\item[(ii)] no supplemented chief factor of $L$ above $C$ or below $R$ is $L$-connected to $A/B$.
\end{itemize}
In other words, there are $r$ ideals $A_1, \ldots , A_r$ of $L$ such that
\[ C/R = A_1/R \oplus \ldots \oplus A_r/R
\]
where $A_i/R$ is a supplemented chief factor of $L$ which is $L$-connected to $A/B$ for $i=1, \ldots ,r$ and $r$ is the number of supplemented chief factors of $L$ which are $L$-connected to $A/B$ in each chief series for $L$. Moreover, $\phi(L/R)=0$.
\end{theorem}

\begin{corollary}\label{c:Lcon} Two supplemented chief factors of $L$ define the same crown if and only if they are $L$-connected.
\end{corollary}

\begin{theorem}\label{t:comp} Let $L$ be a solvable Lie algebra, and let $C/R=\bar{C}$ be the crown associated with a supplemented chief factor of $L$. Then $\bar{C}$ is complemented in $\bar{L}$, and any two complements are conjugate by an automorphism of the form $1+ ad\,a$ for some $a \in \bar{C}$.
\end{theorem}

Finally, in \cite{barnes}, Barnes determined for a solvable Lie algebra which irreducible $L$-modules $A$ have the property that $H^1(L,A)=0$.

\begin{theorem} Let $L$ be a solvable Lie algebra and let $A$ be an irreducible $L$-module. Then $H^1(L,A)=0$ if and only if $L$ has no complemented chief factor isomrphic to $A$.
\end{theorem}

\section{Covering and Avoidance}
 The subalgebra $U$ {\em avoids} the factor algebra $A_i/A_{i-1}$ if $U \cap A_i = U \cap A_{i-1}$; likewise, $U$ {\em covers} $A_i/A_{i-1}$ if $U + A_i = U + A_{i-1}$. We say that $U$ has the covering and avoidance property of $L$ if $U$ either covers or avoids every chief factor of $L$. We also say that $U$ is a $CAP$-subalgebra of $L$. Then these subalgebras give characterisations of solvable and supersolvable Lie algebras; details can be found in \cite{cap}.
\par

There are a number of ways in which $CAP$-subalgebras arise. For a subalgebra $B$ of $L$ we denote by $[B:L]$ the set of all subalgebras $S$ of $L$ with $B \subseteq S \subseteq L$, and by $[B:L]_{max}$ the set of maximal subalgebras in $[B:L]$; that is, the set of maximal subalgebras of $L$ containing $B$. We define the set $\mathcal{I}$ by $i \in \mathcal{I}$ if and only if $A_i/A_{i-1}$ is not a Frattini chief factor of $L$. For each $i \in \mathcal{I}$ put
\[ \mathcal{M}_i = \{ M \in [A_{i-1}, L]_{max} \colon A_i \not \subseteq M\}.
\]
Then $U$ is a {\em prefrattini} subalgebra of $L$ if 
\[ U = \bigcap_{i \in \mathcal{I}} M_i \hbox{ for some } M_i \in \mathcal{M}_i.
\]
It was shown in \cite{prefrat} that, when $L$ is solvable, this definition does not depend on the choice of chief series, and  that the prefrattini subalgebras of $L$ cover the Frattini chief factors and avoid the rest; that is, they are $CAP$-subalgebras of $L$.
\par

Further examples were given by Stitzinger in \cite{stit}, where he proved the following result (see \cite{stit} for definitions of the terminology used).

\begin{theorem}(\cite[Theorem 2]{stit} Let ${\mathcal F}$ be a saturated formation of solvable Lie algebras, and let $U$ be an ${\mathcal F}$-normaliser of $L$. Then $U$ covers every ${\mathcal F}$-central chief factor of $L$ and avoids every ${\mathcal F}$-eccentric chief factor of $L$.
\end{theorem}

The chief factor $A_i/A_{i-1}$ is called {\em central} if $[L,A_i] \subseteq A_{i-1}$ and {\em eccentric} otherwise. A particular case of the above result is the following theorem, due to Hallahan and Overbeck.

\begin{theorem}\label{t:ho} (\cite[Theorem 1]{ho}) Let $L$ be a metanilpotent Lie algebra. Then $C$ is a Cartan subalgebra of $L$ if and only if it covers the central chief factors and avoids the eccentric ones.
\end{theorem}

A subalgebra $U$ of $L$ will be called {\em ideally embedded} in $L$ if $I_L(U)$ contains a Cartan subalgebra of $L$, where $I_L(U) = \{ x \in L : [x,U] \subseteq U \}$ is the {\em idealiser} of $U$ in $L$ . Clearly, any subalgebra containing a Cartan subalgebra of $L$ and any ideal of $L$ is ideally embedded in $L$. Then we have the following extension of Theorem \ref{t:ho}.

\begin{theorem}\label{t:ie} Let $L$ be a metanilpotent Lie algebra and let $U$ be ideally embedded in $L$. Then $U$ is a $CAP$-subalgebra of $L$.
\end{theorem}

\begin{corollary}\label{c:iesolv} Let $L$ be any solvable Lie algebra and let $U$ be an ideally embedded subalgebra of $L$ with $K = N_2(L) \subseteq U$. Then $U$ is a $CAP$-subalgebra of $L$.
\end{corollary}

Another set of examples of $CAP$-subalgebras, which don't require $L$ to be solvable, is given by the next result.

\begin{theorem}Let $L$ be any Lie algebra, let $U$ be a supplement to an ideal $B$ in $L$, and suppose that $B^k \subseteq U$ for some $k \in \mathbb{N}$. Then $U$ is a $CAP$-subalgebra of $L$.
\end{theorem}

We can calculate the dimension of $CAP$-subalgebras in terms of the chief factors that they cover.

\begin{lemma} Let $U$ be a $CAP$-subalgebra of $L$, let $ 0=A_0 < A_1 < \ldots < A_n = L$ be a chief series for $L$ and let ${\mathcal I} = \{i : 1 \leq i \leq n,  U \hbox{ covers } A_i/A_{i-1} \}$. Then $ \dim U = \sum_{i \in {\mathcal I}} (\dim A_i - \dim A_{i-1})$.
\end{lemma}

We have the following characterisations of solvable and supersolvable Lie algebras.

\begin{theorem}  Every one-dimensional subalgebra of $L$ is a $CAP$-subalgebra of $L$ if and only if $L$ is supersolvable.
\end{theorem}

\begin{theorem}\label{t:max}  Let $L$ be a Lie algebra over any field $F$. Then $L$ is solvable if and only if all of its maximal subalgebras are $CAP$-subalgebras.
\end{theorem}

\begin{theorem}\label{t:solv}  Let $L$ be a Lie algebra over a field $F$ which has characteristic zero, or is algebraically closed field and of characteristic greater than $5$. Then $L$ is solvable if and only if there is a maximal subalgebra $M$ of $L$ such that $M$ is a solvable $CAP$-subalgebra of $L$.
\end{theorem}

\end{document}